\documentclass[11pt,twoside,reqno]{amsart}
\usepackage[usenames,dvipsnames]{color}
\usepackage{amsmath,amsthm,amsopn,amssymb,a4wide,bbm,times,enumerate,varioref}
\usepackage[numbers,sort&compress]{natbib}
\usepackage[mathscr]{eucal} 

\numberwithin{equation}{section}

\newtheorem{theorem}{Theorem}[section]
\newtheorem{lemma}[theorem]{Lemma}

\newtheorem{proposition}[theorem]{Proposition}

\newtheoremstyle{note}
     {}
     {}
     {}
     {}
     {\bfseries} 
     {.}
     {.5em}
     {}

\theoremstyle{note}

\newtheorem{remark}[theorem]{Remark}

\def\TT{\Psi}

\def\CT{{C_b(\Psi)}}

\newcommand{\PsF}{\Psi|_F}

\newcommand{\HP}{{H^\infty({\mathcal K}_\Psi)}}

\newcommand{\normbw}[1]{\|#1 \|_{\HP}}

\newcommand{\ip}[2]{{
    \left<
      #1,#2
    \right>}}

\title{Test Functions, Kernels, Realizations and Interpolation}

\author[M.~A.~Dritschel and 
S.~McCullough]{Michael A. Dritschel$^1$
  and Scott McCullough$^2$}

\address{School of Mathematics and Statistics\\
  Merz Court,\\ University of Newcastle upon Tyne\\
  Newcastle upon Tyne\\
  NE1 7RU\\
  UK}

\email{m.a.dritschel@ncl.ac.uk}

\address{Department of Mathematics\\
  University of Florida\\
  Box 118105\\
  Gainesville, FL 32611-8105\\
  USA}

\email{sam@math.ufl.edu}

\subjclass[2000]{47A57 (Primary), 47L55, 47L75, 47D25, 47A13, 47B38,
  46E22 (Secondary)}

\keywords{Agler-Schur class, interpolation, transfer functions,
  multiply connected domains, multivariable, Nevanlinna-Pick,
  Stone-\v{C}ech compactification}

\thanks{${}^1$Research supported by the EPSRC.  \quad ${}^2$Research
  supported by the NSF}

\begin{document}

\begin{abstract}
  Jim Agler revolutionized the area of Pick interpolation with his
  realization theorem for what is now called the Agler-Schur class for
  the unit ball in $\mathbb C^d$.  We discuss an extension of these
  results to algebras of functions arising from test functions and the
  dual notion of a family of reproducing kernels, as well as the
  related interpolation theorem.  When working with test functions,
  one ideally wants to use as small a collection as possible.
  Nevertheless, in some situations infinite sets of test functions are
  unavoidable.  When this is the case, certain topological
  considerations come to the fore.  We illustrate this with examples,
  including the multiplier algebra of an annulus and the infinite
  polydisk.
\end{abstract}

\maketitle

\today

\thispagestyle{empty}

\section{Introduction}
\label{sec:introduction}

Let $\mathbb D^d$ denote the $d$-polydisk in $\mathbb C^d$,
\begin{equation*}
  \mathbb D^d=\left\{z=(z_1,\dots,z_d)\in\mathbb C^d: |z_j|<1, \ 
  j=1,2,\dots,d \right\}.
\end{equation*}
The Agler-Schur class $\mathcal S_d$ consists of those functions
$f:\mathbb D^d\to \mathbb C$ for which there exist positive (that is,
positive semidefinite) kernels $\Gamma_j:\mathbb D^d\times \mathbb D^d
\to \mathbb C$ such that
\begin{equation}
  \label{eq:1}
  1-f(z)f(w)^*=\sum_1^d \Gamma_j(z,w)(1-z_jw_j^* ).
\end{equation}
Here $\zeta^*$ denotes the conjugate of the complex number $\zeta$.

A $\mathbb D^d$ unitary colligation is a pair $\Sigma=(U,\mathcal
E)$, where $\mathcal E=\bigoplus_1^d \mathcal E_j$ is an auxiliary
Hilbert space and
\begin{equation*}
  U=\begin{pmatrix} A &  B\\ C & D\end{pmatrix} :
  \begin{matrix} \mathcal E \\ \oplus \\ \mathbb C \end{matrix}
  \rightarrow 
  \begin{matrix}\mathcal E \\ \oplus \\ \mathbb C \end{matrix}
\end{equation*}
is unitary.  With respect to the decomposition of $\mathcal E$, let
$Z=\bigoplus z_j I_{\mathcal E_j}$, the $d\times d$ block
diagonal operator matrix with diagonal entries $z_jI_{\mathcal E_j}$.
The transfer function $W_\Sigma:\mathbb D^d\to \mathbb C$ of the
colligation $\Sigma$ is defined as
\begin{equation*}
  W_\Sigma = D+ CZ(I-AZ)^{-1}B.
\end{equation*}
 
The \textit{primum mobile} for contemporary work on multivariable
realization with applications to interpolation theory is the following
result of Agler, stated for what is now called the Agler-Schur class
functions on the polydisk.  See \cite{MR1207393}, \cite{MR1665697},
\cite{MR1882259}, \cite{MR1846055}.  This theorem has been generalized in
many directions (\cite{MR1846055} \cite{MR1722812}, \cite{MR1885440},
\cite{MR1797710}, \cite{MR2069781}, to give a few instances).  An operator $C$
on a Hilbert space $H$ is a strict contraction provided $\|C\|<1$.

\begin{theorem}[\cite{MR1207393}]
  \label{thm:getstarted}
  Suppose $f:\mathbb D^d\to \mathbb C$.  The following are equivalent.
  \begin{enumerate}[(i)]
  \item If $k$ is a positive kernel on $\mathbb D^d$ and if
    \begin{equation*}
      (z,w)\mapsto (1- z_jw_j^*)k(z,w)
    \end{equation*}
    is a positive kernel on $\mathbb D^d$ for each $j=1,2,\dots,d$,
    then the kernel
    \begin{equation*}
      (z,w)\mapsto (1-f(z)f(w)^*)k(z,w)
    \end{equation*}
    is also positive;
  \item $f\in\mathcal S_d$;
  \item There is a unitary colligation $\Sigma$ so that
    $f=W_\Sigma$; and
  \item For each tuple $T=(T_1,\dots,T_d)$ of commuting strict
    contractions on a Hilbert space $H$,
    \begin{equation*}
      \|f(T_1,\dots,T_d)\|\le 1.
    \end{equation*}
  \end{enumerate}
\end{theorem}

In newer versions $Z$ is replaced by some matrix function on a domain
in $\mathbb C^d$, the domain being determined by those values at which
the norm of the function is less than $1$ (\cite{MR1885440},
\cite{MR2069781}).  As an example, choosing
\begin{equation*}
  Z=\begin{pmatrix} z_1 & z_2 & \dots & z_d \end{pmatrix}
\end{equation*}
on the unit ball $\mathbb B^d=\{z\in\mathbb C^d: \|z\|<1\}$ in
$\mathbb C^d$ gives rise to the Agler-Schur class of functions
associated with the row contractions and the Nevanlinna-Pick kernel
$k(z,w)=(1-\langle z,w\rangle)^{-1}$.  Of course, choosing $Z$ to be
the diagonal matrix with diagonal entries $z_j$ leads to the
Agler-Schur class of Theorem \ref{thm:getstarted}.

In the paper \cite{DMM} Agler-Schur classes and Agler-Pick
interpolation were considered in the very general setting of algebras
of functions over semigroupoids.  With hopes of mollifying those who
might otherwise be put off by the algebraic formalism of parts of
\cite{DMM}, we drop the semigroupoid structure (or rather, work with a
trivial case).  The realization and interpolation theorems we present
have proofs which in outline follow those found in \cite{DMM}.
However some simplification is achieved in the present setting, and as
a novelty we include a von Neumann type inequality along the lines of
part (iv) of Theorem~\ref{thm:getstarted}.  Furthermore, we show that
strictly contractive functions can be replaced by a certain class
$\mathbb F$ of representations of the algebra generated by the test
functions.  Representations in this family allow for nice
approximations in terms of what we call ``simple representations'' ---
essentially representations involving only finitely many test
functions.  Regarding the interpolation theorem, the realization
formula for the the interpolating function in the Agler-Pick
interpolation theorem has certain additional structure which we
highlight.

As Jim Agler first discovered (\cite{Ag-unpublished}, see also
\cite{MR1882259}), realization and interpolation problems in function
algebras which might not be multiplier algebras for some reproducing
kernel Hilbert space can be effectively handled by means of so-called
test functions (in the realization theorem for the polydisk given
above, these are the coordinate functions).  The test functions are
used to delineate the unit ball in the algebra by means of a duality
with a family of reproducing kernels.  The function algebra is then
the intersection of the multiplier algebras of the reproducing kernel
Hilbert spaces.

The work of Ambrozie \cite{MR2106336} highlighted a second duality at
play in the Agler realization theorem, and in particular in the
definition of the Agler-Schur class.  A similar duality is also
evident in the work of Agler and McCarthy \cite{MR1882259}.  One can view
the test functions as points in some abstract space.  Then in
\eqref{eq:1}, the sum is replaced by a single kernel $\Gamma$
multiplied by $1-E(z)E(w)^*$, where $E(z)$ is the evaluation
functional at the point $z$.  The principle difference between
Ambrozie's approach and that of Agler and McCarthy then comes down to
whether one views the kernel $\Gamma$ as an element of a predual
(Agler and McCarthy), or as an element of a dual space, along with the
introduction of a more general notion of unitary colligation
(Ambrozie).  There is little to distinguish these when there are only
finitely many test functions.  But when there are infinitely many test
functions (something looked at by both), the two approaches are quite
distinct.  We feel that the dual approach offers certain advantages,
which hopefully will be made apparent in what follows.

Our initial motivation was an interest in variants of Theorem~13.8
from \cite{MR1882259} (see also \cite{Ag-unpublished}) which covers
Agler-Pick interpolation on multiply connected domains in $\mathbb C$
(see also \cite{MR532320}, \cite{MR2163865}, \cite{MR1818066},
\cite{MR1386331}, \cite{MR1909298}, and \cite{MR0188824}), and at the
suggestion of John McCarthy \cite{MR1984460}, versions of the Agler
algebra for the infinite polydisk.  While the examples we consider are
replete with analytic structure, a noteworthy feature of the general
theory is that no such structure need be imposed, as will be clear from
the axioms for test functions given below.

Function algebras on multiply connected domains and the infinite
polydisk are examples where infinite families of test functions are
required.  We look at these in some detail, concentrating on the
annulus as the multiply connected domain since all of the salient
features of more general domains are already evident in this example.
In these two cases, the emphasis will be somewhat different.  The
infinite polydisk has a natural choice of test functions.  However
this collection is not compact, the consequence being that there are
functions in the Agler-Schur class that are not simply represented by
some infinite version of \eqref{eq:1}, despite the fact that our
definition of the Agler-Schur class naturally reduces to the original
version for the finite polydisk.  In terms of transfer function
representations for functions in the Agler-Schur class, this is
manifested in the need for the inclusion of representations in the
colligation --- or rather a broader class of representations, since in
fact the decomposition of $\mathcal E$ used in the construction of $Z$
in the transfer function is essentially a representation of a rather
simple form.

The situation is reversed for the annulus in that we are handed a
function algebra (the bounded analytic functions on the annulus), and
the first step then is to find a good set of test functions.  Roughly
speaking, such a set should be as small as possible.  The collection
of test functions we construct in this example is compact, and we show
that it is minimal in the sense that there is no closed subset of this
collection which is also a collection of test functions for this
algebra.

\section{Preliminaries and Main Results}
\label{sec:prel-main-results}

This section contains a discussion of the ingredients that go into the
statement and proof of our generalization of Theorem~\ref{thm:main},
and ends with a statement of the realization formula and Agler-Pick
interpolation theorem.

\subsection{Test Functions and Evaluations}
\label{subsec:test-funct-eval}

For a finite subset $F$ of $X$, let $P(F)$ denote all complex-valued
functions on $X$.  By declaring the indicator functions of points in
$X$ to be an orthonormal basis, $P(F)$ can be identified with the
Hilbert space $\mathbb C^F$.  For a collection $\Psi$ of functions on
$X$, let $\PsF=\{\psi|_F:\psi\in \Psi\} \subset P(F)$.  

A collection $\Psi$ of functions on $X$ is a \textit{collection of
  test functions} provided,
\begin{enumerate}[(i)]
\item For each $x\in X$,
  \begin{equation*}
    \sup\{|\psi(x)|:\psi \in \Psi\} <1; \text{ and}
  \end{equation*}
\item for each finite set $F$, the unital algebra generated by $\PsF$
  is all of $P(F)$.
\end{enumerate}

The second hypothesis, while not essential, simplifies the exposition.
It implies, among other things, that the functions in $\Psi$ separate
the points of $X$.  As we see shortly, the first hypothesis allows us
to use the test functions to define a Banach algebra norm on the
algebra of functions over $X$ with addition and multiplication defined
pointwise.

The function algebras we will be working with may be multiplier
algebras for $H^2(k)$ for some reproducing kernel $k$ (this is what
happens in the classical setting when one studies Nevanlinna-Pick
interpolation), though this is but a special version of what we wish
to consider here.  Rather, test functions will allow us to manage in a
broader context by means of a familiar dual construction.  To this
end, we introduce the algebra of bounded continuous functions over
$\Psi$ with pointwise algebra operations, denoted by $\CT$.

Let $\Psi$ be a collection of test functions.  This is a subspace of
$B(X,\overline{\mathbb D})$, the collection of bounded functions from
$X$ into the closed unit disk $\overline{\mathbb D}$ (equivalently,
$\overline{\mathbb D}^X$), which we endow with the topology of
pointwise convergence.  By Tychonov's theorem, this is a compact
Hausdorff space.  Now $\overline{\mathbb D}$ is a Tychonov space in
the usual metric topology (that is, points are closed and for any
closed set and point disjoint from it, there is a continuous function
separating the two).  Consequently, $\overline{\mathbb D}^X$ is
Tychonov, as is any subspace.  In particular, $\Psi$ is a Tychonov
space.

\textit{A priori} the set $X$ is not assumed to have a topology,
though in fact it inherits one as a subspace of $\CT$.  Since $\Psi$
separates the points of $X$, there is an injective mapping $E:X \to
C_b(\Psi)$, where $E(x) = e_x$ is the evaluation mapping $e_x(\psi) =
\psi(x)$.  Hence we can view $X$ as a subset of the unit ball of
$C_b(\Psi)$.  The same arguments applied above to $\Psi$ now show that
with this topology, $X$ is Tychonov.  In particular, if a net
$\{x_\alpha\}$ converges to $x$ in $X$, then $\psi(x_\alpha)$
converges to $\psi(x)$.  Hence $\Psi$ is a subset of the unit ball of
$C_b(X)$.

By the same token, for a collection of test functions $\Psi$, the
evaluation mapping $F:\Psi\to C_b(\Psi)^*$ can be defined by $F(\psi)
= f_\psi$, where $f_\psi (g) = g(\psi)$ for all $g\in C_b(\Psi)$.  The
mapping $F$ is continuous in the weak-$*$ topology.  Since $\Psi$ is
Tychonov, the map $\Psi \mapsto F(\Psi)$ is a homeomorphism and
$F(\Psi)$ is a subset of the closed unit ball of $C_b(\Psi)^*$, which
is compact.  The Stone-\v{C}ech compactification $\beta\Psi$ of $\Psi$
is then the weak-$*$ closure of $F(\Psi)$.  Since $\beta\Psi$ is a
compact function space, it is pointwise closed, and so contains the
image of the pointwise closure $\overline{\Psi}$ of $\Psi$.  But then
$F$ extends to a homeomorphism of $\overline{\Psi}$ to
$\overline{F(\Psi)}$.  Therefore we can identify $\beta\Psi$ with
$\overline{\Psi}$.

By the way, none of the above depends on the assumption that the test
functions are complex valued, with the exception of the conclusion
that $\CT$ is an algebra generated by $\{E(x): x\in X\}$.  We could,
for example, take our test functions to have values in a Hilbert
$C^*$-module $\mathcal M$ which we may concretely view as a norm
closed subalgebra of $B(\mathcal H, \mathcal K)$, the bounded
operators between Hilbert spaces $\mathcal H$ and $\mathcal K$.  This
is a corner of a $C^*$-subalgebra $\mathcal C$ of $B(\mathcal H \oplus
\mathcal K)$ with the property that each element of $\mathcal M$ is
the corner of some element of $\mathcal C$ with the same norm.  Since
this added generality might obscure the main thrust of the paper, we
restrict our attention to the simpler setting of scalar valued test
functions.

\subsection{Kernel/test function duality and the Agler-Schur
class}
\label{subsec:kern-funct-dual}

Let $\mathcal B$ denote a $C^*$-algebra with Banach space dual
$\mathcal B^*$.  A positive (that is, positive semidefinite) kernel on
a subset $Y$ of $X$ with values in $\mathcal B^*$ is a function
$\Gamma:Y\times Y\to \mathcal B^*$ such that for any finite subset $F$
of $Y$ and function $f:F\to \mathcal B$
\begin{equation*}
  \sum_{a,b\in F} \Gamma(a,b)(f(b)^*f(a)) \ge 0.
\end{equation*}
Naturally the restriction of a positive kernel on $X$ to a subset $F$
is still positive.  We use $M(F,\mathcal B^*)^+$ to denote the
collection of positive definite kernels on $F\subseteq X$ with values
in $\mathcal B^*$.  In the case that $\mathcal B = \mathbb C$ the
modifier \textit{``with values in''} is dropped.

Given a collection $\Psi$ of test functions, write $\mathcal K_\Psi$
for the collection of positive kernels $k$ on $X$ such that for each
$\psi\in \Psi$ the kernel,
\begin{equation*}
  X\times X \ni (x,y) \mapsto (1-\psi(x)\psi(y)^*)k(x,y)
\end{equation*}
is positive.  This is nonempty, since it contains the kernel which is
identically $0$.  Condition (i) in the definition of a collection of
test functions implies that the so-called \textit{Toeplitz kernel} $s$
with $s(x,x) = 1$ for all $x$ and $s(x,y) = 0$ if $y\neq x$ is also in
$\mathcal K_\Psi$.

Say that a function $\varphi:X\to \mathbb C$ is in $\HP$ if there is a
real number $C$ so that, for each $k\in\mathcal K_\Psi$, the kernel
\begin{equation}
  \label{eq:inHP}
  X\times X \ni (x,y) \mapsto (C^2-\varphi(x)\varphi(y)^*)k(x,y)
\end{equation}
is positive, in which case write $C_\varphi$ for the infimum over all
such $C$ (so $C_\varphi$ is independent of $k$).  Then
\begin{equation*}
  \normbw{\varphi} = C_\varphi
\end{equation*}
defines a Banach algebra norm on $\HP$.  Since the Toeplitz kernel $s$
is in $\mathcal K_\Psi$, norm convergent sequences in $\HP$ converge
pointwise, and since positivity of the kernels in \eqref{eq:inHP} is
verified on finite sets, completeness is easily checked.  By
definition, the test functions are in the unit ball of $\HP$.  Indeed,
if $\mathcal K_\Psi$ consists solely of those kernels which are
conjugate equivalent to the Toeplitz kernel (that is, each
$k\in\mathcal K_\Psi$ has the form $k(x,y) = c(x)s(x,y)c(y)^*$ for
some function $c$) then we then have $\HP = C_b(X)$

Replacing $\Psi$ by its pointwise closure (equivalently, $\beta\Psi$)
adds nothing new.  We end up with the same collection of kernels (that
is, $\mathcal K_{\beta\Psi}=\mathcal K_\Psi$), and hence the same
space of functions $\HP$ with the same topology.

A function $f:X \to \mathbb C$ is said to be in the
\textit{Agler-Schur class} if there exists a positive kernel $\Gamma:
X\times X \to \CT^*$ so that for all $x,y \in X$,
\begin{equation*}
  1-f(x)f(y)^*= \Gamma(x,y)(1-E(x)E(y)^*).
\end{equation*}
In the case that $X = \mathbb D^d$ and $\Psi$ consists of the
coordinate functions, we recover the original definition of the
Agler-Schur class from the Introduction.

\subsection{$\CT$-unitary colligations, transfer functions, and the
  class $\mathbb F$}
\label{subsec:ct-unit-coll}

For a collection of test functions $\Psi$, following \cite{MR2106336},
define a \textit{$\CT$-unitary colligation} $\Sigma$ to be a tuple
$\Sigma= (U,\mathcal E,\rho)$, $\mathcal E$ a Hilbert space, $U$
unitary on $\mathcal E \oplus \mathbb C$, and
\begin{equation*}
  \rho:\CT\to  B(\mathcal E)
\end{equation*}
a unital $*$-representation.  Writing $U = \begin{pmatrix} A & B \\ C
  & D \end{pmatrix}$, the \textit{transfer function}
associated to $\Sigma$ is defined as
\begin{equation}
  \label{eq:transfer-bw2}
  W_\Sigma(x)= D+CZ(x)(I-AZ(x))^{-1}B,
\end{equation}
where $Z(x) = \rho(E(x))$.  As observed earlier, by assumption (i) for
test functions, $\|E(x)\| < 1$ for all $x\in X$, and since $\rho$ is
unital, it is contractive.  Hence $\|Z(x)\| < 1$ for all $x\in X$ and
the definition of $W_\Sigma$ makes sense.  Additionally, as the next
lemma indicates, $W_\Sigma$ is contractive.

\begin{lemma}
  \label{lem:tfs-are-contractive}
  Let $\Sigma= (U,\mathcal E,\rho)$ be a $\CT$-unitary colligation
  with associated transfer function $W_\Sigma$.  Then $\|W_\Sigma\|
  \leq 1$.
\end{lemma}

\begin{proof}
  This is a standard calculation.  Using the relations between the
  elements of $U = \begin{pmatrix} A & B \\ C & D \end{pmatrix}$
  implied by the assumption that it is unitary, we find that for any
  $k\in \mathcal K_\Psi$,
  \begin{equation*}
    \begin{split}
      &(1-W_\Sigma(x) W_\Sigma(y)^*)k(x,y) = k(x,y) -
      W_\Sigma(x)k(x,y) W_\Sigma(y)^* \\
       = &C(I-Z(x)A)^{-1}(k(x,y) - Z(x)k(x,y)Z(y)^*)
       (I-Z(y)A)^{-1\,*}C^*.
    \end{split}
  \end{equation*}
  Given a finite set $F$, the matrix
  \begin{equation*}
    (k(x,y) - E(x)k(x,y)E(y)^*)_{x,y\in F}
  \end{equation*}
  is positive, since its value at a test function $\psi$ is
  \begin{equation*}
    (k(x,y) - \psi(x)k(x,y)\psi(y)^*)_{x,y\in F}.
  \end{equation*}
  Since $\rho$ is contractive, the result follows.
\end{proof}

A unital representation $\pi:\HP\to B(\mathcal H)$ is \textit{weakly
  continuous} if whenever $\varphi_\alpha$ is a bounded net from $\HP$
which converges pointwise to a $\varphi\in\HP$, then
$\pi(\varphi_\alpha)$ converges in the weak operator topology to
$\pi(\varphi)$.

We say that $\pi$ is in the class $\mathbb F$ if
$\pi$ is a representation of $H^\infty(\mathcal K_\Psi)$ on a Hilbert
space $\mathcal H$ such that 
\begin{enumerate}[(i)]
\item $\pi$ is weakly continuous; and 
\item $\pi$ is contractive on test functions (that is, $\|\pi(\psi)\|\le 1$
  for each $\psi\in\Psi$).
\end{enumerate}
An example of such a representation is one for which $\|\pi(\psi)\| <
1$ for each $\psi\in\Psi$.  However the class $\mathbb F$ includes
somewhat more, since for example, it also contains the identity
representation, $\pi(\psi) = \psi$ for all $\psi\in\Psi$.
The main fact about this class is that these representations
are automatically contractive and approximately respect
transfer function representations.

\subsection{Abstract realization and Agler-Pick interpolation}
\label{subsec:mainresult}

The following theorem is the analogue of Theorem~\ref{thm:getstarted}
and in fact contains that theorem as a special case when $X = \mathbb
D^d$ and $\Psi$ consists of the $d$ coordinate functions.  More or
less a corollary of this is the Agler-Pick interpolation theorem which
follows it.  With the exception of (iv) in Theorem \ref{thm:main} and
the concrete form of the space $\mathcal E$ for the unitary
colligation in Theorem~\ref{thm:APint}, they are in fact special cases
of results to be found in \cite{DMM}.  Note that a corollary of
Theorem~\ref{thm:main} is that the representations in $\mathbb F$ are
contractive (recall that \textit{a priori} only their behavior on test
functions is prescribed).

\begin{theorem}
  \label{thm:main}
  Suppose $\Psi$ is a collection of test functions.  The following are
  equivalent:
  \begin{itemize}
  \item[(i)] $\varphi\in\HP$ and $\normbw{\varphi}\le 1$;
  \item[(iiF)] For each finite set $F\subset X$ there exists a
    positive kernel $\Gamma:F\times F\to \CT^*$ so that for all
    $x,y\in F$
    \begin{equation*}
      1-\varphi(x)\varphi(y)^*= \Gamma(x,y)(1-E(x)E(y)^*);
    \end{equation*}
  \item[(iiX)] There exists a positive kernel $\Gamma:X\times X\to
    \CT^*$ so that for all $x,y\in X$
    \begin{equation*}
      1-\varphi(x)\varphi(y)^*= \Gamma(x,y)(1-E(x)E(y)^*);
    \end{equation*}
  \item[(iii)] There is a colligation $\Sigma$ so that
    $\varphi=W_\Sigma$;
  \item[(ivS)] For every representation $\pi$ of $\HP$ such that
    $\|\pi(\psi)\| < 1$ for all $\psi\in\Psi$, $\|\pi(\varphi)\| \leq
    1$; and
  \item[(ivF)] For every representation of $\HP$ in $\mathbb F$,
    $\|\pi(\varphi)\| \leq 1$.
  \end{itemize}
\end{theorem}

The Agler-Pick interpolation theorem corresponding to Theorem
\ref{thm:main} is the following.

\begin{theorem}
  \label{thm:APint}
  Suppose $\TT$ is a collection of test functions.  Let $F$, a finite
  subset of $X,$ and $\xi:F\to \mathbb D$ be given.
  \begin{itemize}
  \item[(i)] There exists $\varphi\in\HP$ so that $\normbw{\varphi}\le
    1$ and $\varphi|_F=\xi$;
  \item[(ii)] for each $k\in\mathcal K_\Psi$, the kernel
    \begin{equation*}
      F\times F \ni (x,y) \mapsto (1-\xi(x)\xi(y)^*)k(x,y)
    \end{equation*}
    is positive;
  \item[(iii)] there exists a positive kernel $\Gamma:F\times F\to
    \CT^*$ so that for all $x,y\in F$
    \begin{equation}
      \label{eq:lurk}
      1-\xi(x)\xi(y)^*= \Gamma(x,y)(1-E(x)E(y)^*).
    \end{equation}
  \end{itemize}
  
  Moreover, in this case we also have $\Psi$ is compact (say replacing
  $\Psi$ by its closure) there is a bounded positive measure $\mu$ on
  $\Psi$ so that the interpolant $\varphi$ has a transfer function
  realization of the form
  \begin{equation*}
    \varphi(x)=D+CE(x)(I-AE(x))^{-1}B
  \end{equation*}
  for a unitary
  \begin{equation*}
    U=\begin{pmatrix} A &B\\ C&D \end{pmatrix} : 
    \begin{matrix} \mathbb C^n\otimes L^2(\mu) \\ \oplus \\ \mathbb C 
    \end{matrix}\rightarrow
    \begin{matrix} \mathbb C^n\otimes L^2(\mu) \\ \oplus \\ \mathbb C
    \end{matrix},
  \end{equation*}
  where the $E(x)$ is interpreted as the multiplication operator on
  $\mathbb C^n\otimes L^2(\mu)$ given by $E(x) h\otimes f = h\otimes
  E(x) f$.
\end{theorem}

\begin{remark}\rm (In the last theorem the representation
  $\rho$ in the definition of $\Psi$-unitary colligation is simply a
  multiple of the representation of $C(\Psi)$ as multiplication
  operators on $L^2(\mu)$.  

  Note that in Theorem~\ref{thm:main} it is certainly not the case
  that all positive $\Gamma$ give rise to a $\varphi$.  Indeed,
  Theorem \ref{thm:APint} says a $\Gamma$ corresponding to an
  Agler-Pick interpolation problem can be chosen with additional
  structure.  This will be highlighted in our discussion of the
  annulus.
\end{remark}

As a consequence of the realization theorem, we also get the
following.

\begin{proposition}
  \label{prop:a-s-class-closed}
  The Agler-Schur class is closed in the topology of pointwise
  convergence.  In particular, it contains the closure of the test
  functions.
\end{proposition}

\subsection{Organization}
\label{subsec:organization}

The remainder of the paper is organized as follows.  Section
\ref{sec:preliminaryresults} collects results needed to prove Theorem
\ref{thm:main}, which is then given in Section \ref{subsec:proof}.
Section \ref{subsec:interpolate} contains the proof of
Theorem~\ref{thm:APint}, the basic Agler-Pick interpolation result
companion to Theorem \ref{thm:main}.  Examples are found in Section
\ref{sec:examples}.  These include the annulus algebra and the Agler
algebra of the infinite polydisk.  The case of the annulus $\mathbb A$
is covered in the greatest detail, and it is shown that modulo natural
equivalences, the set of scalar-valued inner functions defined on
$\mathbb A$ with precisely two zeros in $\mathbb A$ and a particular
normalization is a minimal collection of test functions for
$H^\infty(\mathbb A)$.  Further examples show that a certain amount of
care is needed in working in the context presented here, especially
when the collection of test functions is not finite.

\section{Ingredients}
\label{sec:preliminaryresults}

This section collects results preliminary to the proofs of
Theorems~\ref{thm:main} and~\ref{thm:APint}.

\subsection{Simple representations and the class $\mathbb F$}
\label{subsec:reps}

When dealing with an infinite collection of test functions, it is
often useful to approximate using a finite subset.  This is the idea
behind simple representations, which are central to our proof that
(iii) implies (iv) in Theorem~\ref{thm:main}, where we approximate the
function $x\mapsto \rho(E(x))$ appearing in the transfer function
realization.

To be more precise, given a collection of test functions $\Psi$ with
$\psi_j\in \Psi$, $j=1,\ldots, N$, and orthogonal projections $P_j$
such that $\sum_{j=1}^N P_j = I$, define a \textit{simple
  representation} $\rho:\CT \to B(\mathcal E)$ to have the form
\begin{equation*}
  \rho(f) = \sum_{j=1}^N P_j f(\psi_j).
\end{equation*}
Clearly $\rho$ is unital.

Set $Z(x) = \sum_{j=1}^N P_j \psi_j(x) = \rho(E(x))$ and suppose $\pi$
is a representation of $\HP$ on $B(\mathcal H)$.  Then it is natural
to define
\begin{equation}
  \label{eq:defn-piZ}
  \pi(Z) = \sum_{j=1}^N P_j \otimes \pi(\psi_j).
\end{equation}

\begin{lemma}
  \label{lem:simple}
  If
  \begin{enumerate}[(i)]
  \item $\Sigma = (U,\mathcal E,\rho)$ is a $\CT$-unitary
    colligation;
  \item the representation $\rho$ is simple;
  \item $\pi \in \mathbb F$; and
  \item $\varphi= D+CZ(I-A Z)^{-1}B$, where $Z = \rho(E)$,
  \end{enumerate}
  then $\pi(\varphi)$ is a contraction.
\end{lemma}

\begin{proof}
  Let $0<r<1$ and define
\begin{equation*}
  \varphi_r = D+CrZ(I-rA Z)^{-1}B.
\end{equation*}
  
  Fix $x\in X$.  Since $E(x)$ is a strict contraction and $\rho$ is a
  contractive representation, $Z = \rho(E(x))$ is also a strict
  contraction.  It follows that pointwise, $\sum_1^M (rAZ(x))^n$
  converges in norm to $(I-rAZ(x))^{-1}$.  Consequently,
  \begin{equation*}
    \varphi_r^M = D + C Z \sum_1^M  (AZ)^n B,
  \end{equation*}
  converges pointwise with $M$ to $\varphi_r$ and so the sequence
  $\pi(\varphi_r^M)$ converges weakly to $\pi(\varphi_r)$.

  Let $A_j=AP_j$ and for an $n$-tuple $\alpha=(\alpha_1,\dots,\alpha_n)$
  with each $\alpha_j \in \{1,2,\dots,N\}$, let
  $A_\alpha = A_{\alpha_1} A_{\alpha_2} \cdots A_{\alpha_n}$. 
  Define $\psi_\alpha$ similarly.    Let $|\alpha|=n$.
  By expanding $\varphi_r^M$ we have
 \begin{equation*}
    \varphi_M = D + \sum_j\sum_{|\alpha|\le M} r^{|\alpha|+1} CP_j
    A_\alpha B \psi_j \psi_\alpha.
 \end{equation*}
  Thus
 \begin{equation*}
    \begin{split}
      \pi(\varphi_r^M) =& D\otimes I + \sum_j\sum_{|\alpha|\le M}
      r^{|\alpha|+1} CP_j A_\alpha B \pi(\psi_j) \pi(\psi)_\alpha. \\
      =& D\otimes I + (C\otimes I)r\pi(Z) \sum_{n=1}^M ((r A\otimes
      I)\pi(Z))^n (B\otimes I),
    \end{split}
 \end{equation*}
  according to the definition of $\pi(Z)$ in \eqref{eq:defn-piZ}.
  The right side converges in norm with $M$ to
  \begin{equation*}
    D\otimes I +(C\otimes I)r\pi(Z) \left ( (I\otimes I)- r (A\otimes I)
      \pi(Z) \right )^{-1} (B\otimes I),
  \end{equation*}
  giving a transfer function
  representation for $\pi(\varphi_r)$.  The proof that
  $\pi(\varphi)$ is a contraction proceeds along the lines of that
  given for Lemma~\ref{lem:tfs-are-contractive} and makes use
  of the assumption (built into the definition of $\pi$) that
  each $\pi(\psi_j)$ is a contraction.
  
  To complete the proof, note that $\varphi_r$ converges (with $r$)
  pointwise boundedly to $\varphi$ and thus $\pi(\varphi_r)$ converges
  WOT to $\pi(\varphi)$. Since each $\pi(\varphi_r)$ is contractive,
  so is $\pi(\varphi)$.
\end{proof}

\begin{proposition}
   Suppose $\pi\in \mathbb F$.
  \label{prop:pw-simple}
  If
  \begin{enumerate}[(i)]
  \item $\Sigma = (U,\mathcal E,\rho)$ is a $\CT$-unitary colligation;
    and
  \item $\varphi= D+CZ(I-A Z)^{-1}B$, where $Z = \rho(E)$,
  \end{enumerate}
  then there exists a net of simple representations $\rho_\alpha:\CT
  \to B(\mathcal E)$ such that the net
  \begin{equation*}
    \varphi_\alpha = D+CZ_\alpha(I-A Z_\alpha)^{-1}B,
  \end{equation*}
  converges pointwise to $\varphi$, where $Z_\alpha = \rho_\alpha(E)$.
  Consequently, $\pi(\varphi_\alpha)$ converges weakly to
  $\pi(\varphi)$, and so $\pi(\varphi)$ is a contraction.
\end{proposition}

\begin{proof}
  Consider the collection $\mathcal F$ consisting of ordered pairs
  $(F,\epsilon)$ where $F$ is a finite subset of $X$ and $\epsilon >0$
  ordered by $(F,\epsilon)\le (G,\delta)$ if $F\subset G$ and
  $\delta\le \epsilon$.  With this order $\mathcal F$ is a directed
  set.

  Given $\alpha=(F,\epsilon)\in\mathcal F$, by the compactness of
  $\beta\Psi$ there exists a finite collection $\mathcal
  U=\{U_1,\dots, U_m\}$ of nonempty open sets which covers $\beta
  \Psi$ with the property that for any two
  $\psi^\prime,\psi^{\prime\prime} \in U_j$ and any $x\in F$,
  \begin{equation}
    \label{eq:simple-part}
    |\psi^\prime(x)-\psi^{\prime\prime}(x)|<\epsilon.
  \end{equation}

  We construct a partition $\Delta_\alpha=\{\Delta_1,\dots,\Delta_m\}$
  of $\beta \Psi$ from $\mathcal U$ in the usual way.  Let
  \begin{equation*}
    \Delta_m = U_m \setminus (\cup_{j=1}^{m-1} U_j).
  \end{equation*}
  Choose points $\psi_j^\alpha \in U_j \cap \Psi$.  While
  $\psi_j^\alpha$ need not be in $\Delta_j$ it is the case that if
  $\psi\in\Delta_j$, then \eqref{eq:simple-part} holds with
  $\psi^\prime =\psi_j$ and $\psi^{\prime\prime} = \psi$.  Since, for
  $x\in F$,
  \begin{equation}
    \label{eq:2}
    \left\| \sum \psi_j^\alpha(x) Q(\Delta_j) -\int \psi(x)dQ(\psi)
    \right\| \le \sum \left\| \int_{\Delta_j}
    (\psi_j^\alpha(x)-\psi(x)) dQ(\psi) \right\| \le \epsilon 
  \end{equation}

  Let $Q$ denote spectral measure associated to $\rho$ so that
  \begin{equation*}
    \rho(f)=\int_{\beta\Psi} f\,dQ.
  \end{equation*}
  Define
  \begin{equation*}
    \rho_\alpha (f) =\sum_{j=1}^n f(\psi_j)Q(\Delta_j)
  \end{equation*}
  and let $Z_\alpha (x)=\rho`_\alpha (E(x))$.  It follows by
  \eqref{eq:2} that for $\alpha = (F,\epsilon)$, $\|Z_\alpha(x)
  -Z(x)\| \le \epsilon$ for $x\in F$ and this remains true if $\alpha$
  is replaced by any $\beta \ge \alpha$.  Since $\|E(x)\| < 1$, we
  have $0< \delta = \sup_{x\in F}(1-\|E(x)\|)/2$, and so for $r =
  1-\delta/2$ and $\epsilon < \delta/2$, it follows that
  $\|Z_\alpha(x)\| < r = 1-\epsilon$.

  By Lemma~\ref{lem:tfs-are-contractive}, for any $\alpha$,
  \begin{equation*}
    \varphi_\alpha = D+C Z_\alpha (I- A Z_\alpha)^{-1} B.
  \end{equation*}
  is a contraction.  Note that
  \begin{equation*}
    \begin{split}
      &(AZ_\alpha)^n - (AZ)^n \\ = & A(Z_\alpha - Z) AZ_\alpha \cdots
      AZ_\alpha + AZ A(Z_\alpha - Z) AZ_\alpha \cdots AZ_\alpha 
      + \ldots +
      AZ\cdots AZ A(Z_\alpha - Z),
    \end{split}
  \end{equation*}
  and so
  \begin{equation*}
    \|(AZ_\alpha)^n - (AZ)^n\| \leq n \|A\| \epsilon r^{n-1}.
  \end{equation*}
  Hence for suitably chosen $\alpha$,
  \begin{equation*}
    \begin{split}
      &\|Z_\alpha (I- A Z_\alpha)^{-1} - Z_ (I- A Z)^{-1} \| \\ = &
      \|(Z_\alpha - Z)(I- A Z_\alpha)^{-1} + Z[ (I- A Z_\alpha)^{-1} -
      Z_ (I- A Z)^{-1}] \| \\
      \leq & \epsilon \left[\frac{1}{1-r} + \frac{r^2}{(1-r)^2}\right],
    \end{split}
  \end{equation*}
  Thus the bounded net $\varphi_\alpha$ converges pointwise to
  $\varphi$.

  As constructed each $\varphi_\alpha$ has a simple transfer function
  representation, and so by Lemma \ref{lem:simple},
  $\pi(\varphi_\alpha)$ is a contraction.  Since the net
  $\varphi_\alpha$ is bounded and converges pointwise to $\varphi$,
  then net $\pi(\varphi_\alpha)$ converges in the weak operator
  topology to $\pi(\varphi)$.  Hence $\pi(\varphi)$ is a contraction.
\end{proof}

\subsection{Factorization}
\label{subsec:factorization}

The engine powering the lurking isometry argument in the proof of
(iiX) implies (iii) of Theorem \ref{thm:main} is the factorization in
the following proposition.  A similar result may be found in
\cite{MR2106336}.  A detailed proof, which we hint at, is given in
\cite{DMM}.  See also the book \cite{MR1882259} Theorem 2.53 proof 1.

\begin{proposition}
  \label{prop:factorization}
  If $\Gamma:X\times X \to \CT^*$ is positive, then there exists a
  Hilbert space $\mathcal E$ and a function $L:X\to B(\CT, \mathcal
  E)$ such that
  \begin{equation*}
    \Gamma(x,y)(fg^*)=\ip{L(x)f}{L(y)g}
  \end{equation*}
  for all $f,g\in \CT$.

  Further, there exists a unital $*$-representation $\rho:\CT\to
  B(\mathcal E)$ such that $L(x)ab=\rho(a)L(x)b$ for all $x\in X$,
  $a,b\in\CT$.
\end{proposition}

\begin{proof}[Sketch of the proof]
  Let $V$ denote the vector space with basis $X$.  On the vector space
  $V\otimes \CT$ introduce the positive semidefinite sesquilinear form
  induced by
  \begin{equation*}
    \ip{x \otimes f}{y\otimes g} =\Gamma(x,y)(g^* f),
  \end{equation*}
  where $x,y\in X$ and $f,g\in\CT$.  This is positive semidefinite by
  the hypothesis that $\Gamma$ is positive semidefinite.  Mod out by
  the kernel and complete to get the Hilbert space $\mathcal E$.
  Define $L$ by $L(x)a= x\otimes a$ and verify that this is a bounded
  operator.

  The $*$-representation is induced by the left regular representation
  of $\CT$, $\rho:\CT \to B(\mathcal E)$ with $\rho(a)(x\otimes f)=
  x\otimes af$.  Then check that $\rho$ is a contractive unital
  representation of $\CT$ satisfying $L(x)ab=\rho(a)L(x)b$ for all
  $x\in X$, $a,b\in\CT$.
\end{proof}

\subsection{A closed cone}
\label{subsec:closed-cone}

The proof of (i) implies (ii) in Theorem \ref{thm:main} is based on a
cone separation argument which, in order to work, requires that the
cone be closed and have nonempty interior.  We present some of the
background material here.

Given a finite subset $F\subset X$, let $\CT_F^+$ denote the
collection of positive kernels $\Gamma:F\times F\to \CT^*$.  If
$\Gamma\in\CT_F^+$ and $x\in F$, then $\Gamma(x,x)$ is a positive
linear functional on the unital $C^*$-algebra $\CT$ and therefore,
$\|\Gamma\| = \Gamma(1)$.

Details of the proof outlined below can be found in~\cite{DMM}.

\begin{lemma}
  \label{lem:closedcone}
  Let $\TT$ be a set of test functions.  If for each $x\in X$,
  $\|E(x)\|_\infty <1$, then
  \begin{equation}
    \label{eq:C_F}
    \mathcal C_F=\{ \begin{pmatrix} \Gamma(x,y)(I-E(x)E(y)^*)
    \end{pmatrix}_{x,y\in F}: \Gamma \in \CT_F^+\}
  \end{equation}
  is a closed cone of $|F|\times |F|$ matrices $($where $|F|$ is the
  cardinality of $F)$.
\end{lemma}

\begin{proof}[Sketch of the proof]
  Suppose $M = (\Gamma(x,y)(I-E(x)E(y)^*)) \in \mathcal C_F$.  Since
  $\|E(x)\| < 1$, $1-E(x)E(x)^* > \epsilon 1$ for some $\epsilon > 0$.
  Hence $\frac{1}{\epsilon}M(x,x) \geq \Gamma(x,x) 1 =
  \|\Gamma(x,x)\|$, and so $\|\Gamma(x,x)\| \leq
  \frac{1}{\epsilon}\|M\|$.  Finiteness of $F$ means that there is in
  fact a single $\epsilon$ which will do for all $x\in F$, while
  positivity of $\Gamma$ implies that for $g\in\CT$,
  \begin{equation*}
    2|\Gamma(x,y) g| \leq \Gamma(x,x) 1 + \Gamma gg^* \leq
    \|\Gamma(x,x)\| + \|\Gamma(y,y)\|\|gg^*\|.
  \end{equation*}
  Consequently, $\|\Gamma(x,y)\| \leq \frac{1}{\epsilon}\|M\|$ for all
  $x,y\in F$.

  Now suppose $M_j\in\mathcal C_F$ is a Cauchy sequence.  For each $j$
  there exists $\Gamma_j\in \CT_F^+$ so that
  \begin{equation*}
    M_j=\begin{pmatrix} \Gamma_j(x,y)(I-E(x)E(y)^*)
    \end{pmatrix}_{x,y\in F}.
  \end{equation*}
  Since the $M_j$'s are uniformly bounded, $\Gamma_j(x,y)$ is
  uniformly bounded for all $x,y$ and $j$.  Thus there is a
  subsequence $\Gamma_{j_\ell}$ such that $\Gamma_{j_\ell}(x,y)$
  converges weak-$*$ to $\Gamma(x,y)$.  Likewise,
  $\Gamma_{j_\ell}(x,y)E(x)E(y)^*$ converges weak-$*$ to
  $\Gamma(x,y)E(x)E(y)^*$.  Hence $M=\lim_j M_j =
  (\Gamma(x,y)(1-E(x)E(y)^*)_{x,y\in F}$.  Positivity of $\Gamma$ is a
  consequence of the positivity of the $\Gamma_j$'s.  We conclude that
  $\mathcal C_F$ is closed.
\end{proof}

The next lemma gives an example of a positive kernel in $\CT_F^+$
which will be particularly useful in showing that the cone $\mathcal
C_F$ in \eqref{eq:C_F} has nonempty interior in the subsequent lemma.

\begin{lemma}
  \label{lem:evalmu}
  Let $\TT$ be a set of test functions for which $\|E(x)\|<1$ for all
  $x$.  For each $\psi\in\TT$ the function $\Gamma_\psi:X\times X\to
  \CT^*$ given by
  \begin{equation*}
    \Gamma_\psi(x,y)(f)=\frac{f(\psi)}{1-\psi(x)\psi(y)^*}
  \end{equation*}
  is a positive kernel.  Here $f\in \CT$.
\end{lemma}

\begin{proof}
  Note that $|\psi(x)|\le \|E(x)\|<1$ so that the formula makes sense
  and moreover,
  \begin{equation*}
    S(x,y)=\frac{1}{1-\psi(x)\psi(y)^*}
  \end{equation*}
  defines a positive kernel on $X$.

  With $y=x$,
  \begin{equation*}
    |\Gamma_\psi(x,x)(f)| 
    =\frac{|f(\psi)|}{|1-|\psi(x)|^2} \le \|f\| \frac{1}{|1-|\psi(x)|^2}
  \end{equation*}
  so that $\Gamma_\psi(x,x)$ is indeed in $\CT^*$.

  For a finite set $F\subset X$ and function $f:F\to \CT$
  \begin{equation*}
    \sum_{x,y\in F} \Gamma_\psi(x,y)(f(x)f(y)^*)
    =\sum_{x,y\in F} f(x)(\psi) f(y)(\psi)^* S(x,y) \ge 0,
  \end{equation*}
  since $S$ is a positive kernel on $X$.  It follows that each
  $\Gamma_\psi(x,y)\in\CT^*$ and $\Gamma_\psi$ is a positive kernel on
  $X$.
\end{proof}

\begin{lemma}
  \label{lem:Kpos}
  Let $\TT$ be a set of test functions, $F\subset X$ a finite set, and
  $\mathcal C_F$ the cone in \eqref{eq:C_F}.  Then $\mathcal C_F$
  contains all positive $|F|\times |F|$ matrices, and hence has
  nonempty interior.
\end{lemma}

\begin{proof}
  Let $\Gamma_\psi$ denote the positive kernel from
  Lemma~\ref{lem:evalmu}.  Then
  \begin{equation*}
    [1]=\Gamma_\psi(x,y)(I-E(x)E(y)^*) \in\mathcal C_F,
  \end{equation*}
  where $[1]$ is the matrix with all entries equal to $1$.  For $P$ be
  a positive $|F|\times |F|$ matrix, $\tilde\Gamma$ defined by
  $\tilde\Gamma(x,y) = P(x,y) \Gamma_\psi(x,y)$ is a positive kernel,
  and so $P(x,y) = \tilde\Gamma(x,y)(I-E(x)E(y)^*)$.
\end{proof}

\begin{lemma}
  \label{lem:conjcone}
  The cone $\mathcal C_F$ in \eqref{eq:C_F} is closed under
  conjugation by diagonal matrices; i.e., if $M=(M(x,y))\in\mathcal
  C_F$ and $c:F\to \mathbb C$, then $cMc^*=(c(x)M(x,y)c(y)^*)
  \in\mathcal C_F$.
\end{lemma}

\begin{proof}
  Simply note that if $\Gamma:F\times F\to \CT^*$ is positive, then so
  is $c\Gamma c^*$ defined by $(c\Gamma c^*)(x,y)=c(x)c(y)^*
  \Gamma(x,y)$.
\end{proof}

The next proposition connects the closed unit ball of of $\HP$ with
the cone $\mathcal C_F$.  Further details of the proof sketched below
can be found in Lemmas~5.5 and 3.4 of \cite{DMM}.

\begin{proposition}
  \label{prop:separation}
  Let $F\subset X$ be finite and $\varphi\in\HP$.  If $M_\varphi$
  defined by
  \begin{equation*}
    M_\varphi(x,y) = 1 - \varphi(x) \varphi(y)^*,\qquad x,y\in F,
  \end{equation*}
  is not in $\mathcal C_F$, then there is a kernel $k\in \mathcal
  K_\Psi$ such that the matrix
  \begin{equation*}
    (1 - \varphi(x) \varphi(y)^*) k(x,y))_{x,y\in F}
  \end{equation*}
  is not positive.  That is, $\|\varphi\|_\HP > 1$.
\end{proposition}

\begin{proof}[Sketch of the proof]
  Use a version of the Hahn-Banach theorem to find a linear functional
  $\lambda\neq 0$ on the selfadjoint $|F|\times |F|$ matrices such
  that $\lambda(M) \geq 0$ for all $M\in\mathcal C_F$ but
  $\lambda(M_\varphi) < 0$.  We can find such a $\lambda$ since
  $\mathcal C_F$ is closed and has nonempty interior by
  Lemmas~\ref{lem:closedcone} and~\ref{lem:Kpos}.

  For $f,g\in P(F)$ (viewed as vectors in $\mathbb C^F$), define
  $\ip{f}{g} = \lambda(fg^*)$.  Since $\mathcal C_F$ contains all
  positive $|F|\times |F|$ matrices, this is positive.  Mod out by the
  kernel and call the resulting space $\mathcal H$.  Let $q$ be the
  quotient map.  Show that $\lambda(M_\varphi) < 0$ implies
  $\lambda([1])> 0$, and hence $q(\delta_F) \neq 0$, where $\delta_F$
  is the function in $P(F)$ which is identically $1$.

  Let $\mu$ be a representation of $P(F)$ given by $\mu(g)q(f) =
  q(fg)$, where the product $fg$ is defined pointwise.  Verify that
  $\mu$ is contractive on test functions and that $\mu([1]-\varphi|F
  \varphi^*|F) < 0$.

  What is more, $\delta_F$ is a cyclic vector for $\mu$.  Hence if
  $\xi_x\in P(F)$ is defined to be $1$ at $x$ and zero elsewhere, then
  $\{\ell_x = \mu(\xi_x)\delta_F\}_{x\in F}$ is a basis for $\mathcal
  H$.  Let $\{k_x\}$ be the dual basis, $k(x,y) = \ip{k_x}{k_y}$.
  Then for $c\in \mathbb C$,
  \begin{equation*}
    \begin{split}
      \ip{\mu(c\xi_x)^* k_a}{\ell_b} & = c^*\ip{k_a}{\mu(\xi_x)\ell_b}
      \\
      &= c^*\ip{k_a}{\mu(\xi_x)\mu(\xi_b)\delta_F} \\
      &=
      \begin{cases}
        c^* & \text{if }x=b=a\\
        0 & \text{otherwise}
      \end{cases}
      \\
      &= \ip{c^* k_a}{\ell_b}.
    \end{split}
  \end{equation*}
  So for $f\in P(F)$, $\mu(f)^* k_a = f(a)^* k_a$.  If $f$ is the test
  function $\psi$, this yields that the matrix
  \begin{equation*}
    ((1-\psi(x)\psi(y)^*) k(x,y))_{x,y\in F}
  \end{equation*}
  is positive, while with $f=\varphi$, it is strictly negative.
  Extend $k$ to all of $X\times X$ by setting $k(x,y) = 0$ if either
  $x$ or $y$ are not in $F$.  We then have $((1-f(x)f(y)^*)
  k(x,y))_{x,y\in X}$ is positive when $f$ is a test function (so that
  $k\in\mathcal K_\Psi$), but not positive for $f=\varphi$.
\end{proof}

\subsection{A compact set}
\label{subsec:compactset}

The proof of (iiF) implies (iiX) in Theorem \ref{thm:main} uses
Kurosh's theorem (\cite{MR1882259}, Theorem~2.56), the application of
which requires that certain sets be compact.

Fix $\varphi:X\to \mathbb C$ and a collection of test functions $\TT$.
For $F\subset X$, let
\begin{equation*}
  \Phi_F = \{\Gamma\in\CT_F^+ : 1-\varphi(x)\varphi(y)^*=
  \Gamma(x,y)(1-E(x)E(y)^*) \text{ for } x,y\in F\}.
\end{equation*}
The set $\Phi_F$ is naturally identified with a subset of the
product of $\CT^*$ with itself $|F|^2$ times.

\begin{lemma}
  \label{lem:compact-set}
  If for each $x\in X$, $\|E(x)\|<1$, then the set $\Phi_F$ is
  compact.
\end{lemma}

\begin{proof}
  Let $\Gamma_\alpha$ be a net in $\Phi_F$.  Arguing as in the
  proof of Lemma \ref{lem:closedcone}, we find each
  $\Gamma_\alpha(x,x)$ is a bounded net and thus each
  $\Gamma_{\alpha}(x,y)$ is also a bounded net.  By weak-$*$
  compactness of the unit ball in $\CT^*$ there exists a $\Gamma$ and
  subnet $\Gamma_\beta$ of $\Gamma_\alpha$ so that for each $x,y\in
  F$, $\Gamma_{\beta}(x,y)$ converges to $\Gamma(x,y)$.
\end{proof}

\section{Proofs}
\label{sec:proofs}

We are now set to prove the theorems stated in
Subsection~\ref{subsec:mainresult}.

\subsection{Proof of Theorem \ref{thm:main}}
\label{subsec:proof}

\subsubsection{Proof of (i) implies (iiF)}

Let $\varphi\in\HP$.  If we suppose (iiF) does not hold, then by
Proposition~\ref{prop:separation}, $\|\varphi\|_\HP > 1$.

\subsubsection{Proof of (iiF) implies (iiX)}

The proof here uses Kurosh's Theorem and in much the same way as in
\cite{MR1882259}.

The hypothesis is that for every finite subset $F\subset X$,
$\Phi_F$, as defined in subsection \ref{subsec:compactset} is not
empty, and so by Lemma~\ref{lem:compact-set}, $\Phi_F$ is compact.
For finite set $F\subset G$, define $ \pi_F^G :\Phi_F \to
\Phi_G$ by
\begin{equation*}
  \pi_F^G (\Gamma)  =\Gamma|_{F\times F}.
\end{equation*}
Thus, with $\mathcal F$ equal to the collection of all finite subsets
of $X$ partially ordered by inclusion, the triple
$(\Phi_F,\pi_F^G,\mathcal F)$ is an inverse limit of nonempty
compact spaces.  Consequently, by Kurosh's Theorem, for each
$F\in\mathcal F$ there is a $\Gamma_F \in \Phi_F$ so that whenever
$F,G\in\mathcal F$ and $F\subset G$,
\begin{equation}
  \label{eq:consistent}
  \pi_F^G (\Gamma_G)  =\Gamma_F.
\end{equation}

Define $\Gamma:X\times X\to \CT^*$ by $\Gamma(x,y)=\Gamma_{F}(x,y)$
where $F\in\mathcal F$ and $x,y\in F$.  This is well defined by the
relation in equation (\ref{eq:consistent}).  If $F$ is any finite set
and $f:F\to \CT$ is any function, then
\begin{equation*}
  \sum_{x,y\in F} \Gamma(x,y)(f(x)f(y)^*) =\sum_{x,y\in
    F}\Gamma_F(x,y)(f(x)f(y)^*) \ge 0
\end{equation*}
since $\Gamma_F\in \CT_F^+$.  Hence $\Gamma$ is positive.

\subsubsection{Proof of (iiX) implies (iii)}

Let $\Gamma$ denote the positive kernel of the hypothesis of (iiX).
Apply Lemma \ref{prop:factorization} to find $\mathcal E$, $L:X\to
B(\CT, \mathcal E)$, and $\rho:\CT\to B(\mathcal E)$ as in the
conclusion of the lemma.
  
Rewrite condition (iiX) as
\begin{equation*}
  1+\ip{Z(x)L(x)1}{Z(x)L(x)1} =\varphi(x) \varphi(y)^* +
  \ip{L(x)1}{L(x)1},
\end{equation*}
where we use Proposition~\ref{prop:factorization} to express $L(x)E(x)
= Z(x)L(x)1$ with $Z(x) = \rho(E(x))$.  From here the remainder of the
proof is the standard lurking isometry argument.

Let $\mathcal E_d$ denote finite linear combinations of
\begin{equation*}
  \begin{pmatrix} Z(x)L(x)1 \\ 1 \end{pmatrix} \in \begin{matrix}
    \mathcal E \\ \oplus \\ \mathbb C \end{matrix}
\end{equation*}
and let $\mathcal E_r$ denote finite linear combinations of
\begin{equation*}
  \begin{pmatrix} L(x)1 \\ \varphi(x) \end{pmatrix} \in \begin{matrix}
    \mathcal E \\ \oplus \\ \mathbb C \end{matrix}.
\end{equation*}

Define $V:\mathcal E_d\to \mathcal E_r$ by
\begin{equation*}
  V\begin{pmatrix} Z(x)L(x)1 \\ 1 \end{pmatrix} = \begin{pmatrix}
    L(x)1 \\ \varphi(x) \end{pmatrix},
\end{equation*}
extend by linearity and show that $V$ is a well defined isometry on
$\mathcal E_d$, and hence on $\overline{\mathcal E_d}$.  This further
extends to a unitary operator
\begin{equation*}
  U=\begin{pmatrix} A & B\\ C & D \end{pmatrix} 
  : \begin{matrix} \mathcal H \\ \oplus \\ \mathbb C \end{matrix} \to
  \begin{matrix} \mathcal H \\ \oplus \\ \mathbb C \end{matrix},
\end{equation*}
with $U$ restricted to $\mathcal E_d$ equal to $V$; that is
$U\gamma=V\gamma$ for $\gamma\in\mathcal E_d$.

This gives the system of equations
\begin{equation*}
  \begin{split}
    A Z(x)L(x)1 + B &= L(x)1 \\
    \nonumber C Z(x)L(x)1 + D &= \varphi(x),
  \end{split}
\end{equation*}
which, when solved for $\varphi$, yields
\begin{equation*}
  \varphi(x)=D+CZ(x)(I-AZ(x))^{-1}B,
\end{equation*}
as desired.

\subsubsection{Proof of (iii) implies (ivF)}

This is a direct consequence of Proposition~\ref{prop:pw-simple}.

\subsubsection{Proof that (ivF) is equivalent to (ivS)}

This is trivial in one direction.  In the other, it follows from the
proof of Lemma~\ref{lem:simple}.

\subsubsection{Proof of (ivF) implies (i)}

Take $\pi$ to be the identity representation in
Proposition~\ref{prop:pw-simple}.

\subsection{Agler-Pick Interpolation: Proof of Theorem \ref{thm:APint}}
\label{subsec:interpolate}

It turns out that in the Agler-Pick interpolation setting more can be
said about the transfer function realization of the interpolant.
Suppose $\mu$ is a (positive) measure on $\Psi$.  The functions $E(x)$
determine multiplication operators on $L^2(\mu)$ by the formula
$(E(x)f)(\psi)=\psi(x)f(\psi)$.  Abusing notation, for a positive
integer $n$, let $E(x)$ also denote the operator $I_n\otimes E(x)$ on
$\mathbb C^n \otimes L^2(\mu)$, or more precisely, the representation
$\rho(E(x))=I_n\otimes E(x)$

\begin{proof}[Proof of Theorem~\ref{thm:APint}]
  If $\varphi$ exists, the implication (i) implies (iiF) of Theorem
  \ref{thm:main} applied to $F$ establishes the existence of $\Gamma$.

  Conversely, suppose a positive $\Gamma$ satisfying equation
  \eqref{eq:lurk} exists.  View $\Gamma$ as an $n\times n$ matrix
  \begin{equation*}
    \Gamma =\begin{pmatrix} \Gamma(x_\ell,x_j)\end{pmatrix}_{j,\ell}
  \end{equation*} 
  with entries from $\CT^*$.  The converse can be proved using the
  factorization from Proposition \ref{prop:factorization}.  However
  the proof of the last part about the measure $\mu$ requires a
  somewhat more concrete factorization of $\Gamma$.

  Choose a positive measure $\mu$ on $\TT$ so that each
  $\Gamma(x_\ell,x_j)$ is absolutely continuous with respect to $\mu$.
  We can without loss of generality assume that the measure is defined
  on the closure of $\TT$ (if this is not already closed), and hence
  we may assume that the measure $\mu$ is bounded.  By Radon-Nikodym,
  there exist $L^\infty(\mu)$ functions $F_{j,\ell}$ so that
  $\Gamma(x_\ell,x_j)=F_{\ell,j}\, d\mu$.  In particular, the
  matrix-valued function $F$ can be identified with an element of the
  $C^*$-algebra of $n\times n$ matrices with entries from
  $L^\infty(\mu)$.  The fact that $\Gamma$ is positive implies that
  $F$ is (almost everywhere $\mu$) pointwise positive.  Consequently,
  there exists vectors $H$ from $C^n\otimes L^\infty(\mu)$ so that
  $F_{\ell,j}=H(x_\ell)H(x_j)^*$.  This gives the factorization,
  \begin{equation*}
    \Gamma =HH^* \, d\mu.
  \end{equation*}

  Observe
  \begin{equation*}
    \begin{split}
      \Gamma(x_\ell,x_j)(1-&E(x_\ell)E(x_j)^*)\\
      =& \int H(x_\ell) H(x_j)^*\, d\mu
      -\int H(x_\ell) E(x_\ell)^*E(x_j) H(x_j)^* \, d\mu \\
      =& \langle H(x_\ell),H(x_j) \rangle -\langle E(x_j)H(x_j),
      E(x_\ell)H(x_\ell)\rangle.
    \end{split}
  \end{equation*}
  Thus, equation (\ref{eq:lurk}) becomes,
  \begin{equation}
    \label{eq:lurkspecial}
    1+ \ip{E(x_j)H(x_j)}{E(x_\ell)H(x_\ell)} =
    \xi(x_j)\xi(x_\ell)^* + \ip{H(x_\ell)}{H(x_j)}.
  \end{equation}

  A lurking isometry argument as in the proof of (iiX) implies (iii)
  for Theorem~\ref{thm:main} allows us to define a unitary operator
  \begin{equation*}
    U=\begin{pmatrix} A &B\\ C&D \end{pmatrix} : 
    \begin{matrix}  \mathbb C^n\otimes L^2(\mu) \\ \oplus \\ \mathbb C
    \end{matrix}\rightarrow 
    \begin{matrix}  \mathbb C^n\otimes L^2(\mu) \\ \oplus \\ \mathbb C
    \end{matrix}
  \end{equation*}
  with
  \begin{equation*}
    U
    \begin{pmatrix}
      E(x_j)H(x_j) \\ 1
    \end{pmatrix}
    =
    \begin{pmatrix}
      H(x_j) \\ \xi(x_j)
    \end{pmatrix},
  \end{equation*}
  which can then be solved to give
  \begin{equation*}    
    \xi(x_j)= D+CE(x_j)(I-AE(x_j))^{-1}B.
  \end{equation*}
  Define
  \begin{equation*}
    \varphi(x)=D+CE(x)(I-A E(x))^{-1}B
  \end{equation*}
  for $x\in X$.  Then $\varphi$ extends $\xi$ and the implication
  (iii) implies (i) of Theorem \ref{thm:main} completes the proof.
\end{proof}

\subsection{Proof of Proposition~\ref{prop:a-s-class-closed}}
\label{subsec:closed-class}

  Let $\Psi$ be a collection of test functions, $\mathcal K_\Psi$ and
  $\HP$ as above, and suppose $\varphi_\alpha$ is a net in the
  Agler-Schur class of $\HP$.  Then for all $\alpha$,
  $\|\varphi_\alpha\| \leq 1$, and so $\|\varphi\| \leq 1$.  Fix
  $F\subset X$ finite.  Then there is a $\Gamma_{F,\alpha} \geq 0$
  such that
  \begin{equation*}
    1-\varphi_\alpha(x)\varphi_\alpha(y)^* =
    \Gamma_{F,\alpha}(x,y)(1-E(x)E(y)^*), \qquad x,y\in F.
  \end{equation*}
  So the matrix $M_\alpha = (1-\varphi_\alpha(x)\varphi_\alpha(y)^*)
  \in \mathcal C_F$.  Since by Lemma~\ref{lem:closedcone} $\mathcal
  C_F$ is closed, arguing as at the end of the proof of that lemma, we
  have a $\Gamma_F \geq 0$ such that
  \begin{equation*}
    1-\varphi(x)\varphi(y)^* = \Gamma_F (x,y)(1-E(x)E(y)^*), \qquad
    x,y\in F.
  \end{equation*}
  Applying (iiF) implies (iiX) of Theorem~\ref{thm:main}, it follows
  that $\varphi$ is in the Agler-Schur class of $\HP$.  Finally, since
  the test functions are a subset of the Agler-Schur class, the last
  statement is obvious.

\section{Examples}
\label{sec:examples}

In this section we concentrate on two main examples where an infinite
collection of test functions is required; the annulus and the infinite
polydisk.  We then close with a few further examples illustrating the
necessity of various parts of our definitions of test functions.

\subsection{The annulus}
\label{subsec:annulus}

Fix $q\in (0,1)$ and write $\mathbb A=\mathbb A_q$ for the annulus
\begin{equation*}
  \mathbb A=\{z\in\mathbb C: q<|z|<1\}.
\end{equation*}
Let $H^\infty(\mathbb A)$ denote the bounded analytic functions on
$\mathbb A$.  There is a collection of functions $\vartheta_t$
naturally parameterized by $t$ in the unit circle $\mathbb T$ with the
property that each $\vartheta_t$ is unimodular on the boundary of
$\mathbb A$ (and so extending analytically across the boundary) and
has precisely two zeros in $\mathbb A$.  Moreover, any function with
these properties is, up to pre-composition with an automorphism of
$\mathbb A$ and post-composition with an automorphism of $\mathbb D$,
one of these $\vartheta_t$.  We begin by constructing $\vartheta_t$
and showing that $\Theta = \{\vartheta_t: t\in\mathbb T\}$ is indeed a
family of test functions for $H^\infty(\mathbb A)$.

Let $B_0=\{|z|=1\}$ and $B_1=\{|z|=q\}$ denote the boundary components
of the boundary $B$ of $\mathbb A$.  For normalization, fix a base
point $b\in\mathbb A$ such that $|b|\ne \sqrt{q}$.  Using Green's
functions (or otherwise), for each point $\alpha\in B$ there exists a
unique positive harmonic function $h_\alpha$ whose boundary values
come from the measure on $B$ with point mass at $\alpha$.  If $h$ is
any positive harmonic function on $\mathbb A$ there is a (positive)
measure $\mu$ on $B$ so that
\begin{equation}
  \label{eq:harmonic}
  h(z)=\int_B h_\gamma \, d\mu(\gamma) 
  =\int_{B_0} h_\alpha \, d\mu(\alpha) + 
  \int_{B_1} h_\beta \, d\mu(\beta).
\end{equation}
The harmonic function $h$ is the real part of an analytic function $f$
if and only if $\mu(B_0)=\mu(B_1)$.  In particular, given $\alpha\in
B_0$ and $\beta \in B_1,$ the function $h_\alpha+h_\beta$ is
the real part of an analytic function $f_{\alpha,\beta}$
which we may normalize by requiring the imaginary part of
$f_{\alpha,\beta}(b)=0$.  Since both $h_\alpha$ and $h_\beta$
are nonnegative on the boundary they are both positive inside
the annulus and so $f(b)>0$.
Then because the boundary values for the $h_\alpha$'s are point
masses, \eqref{eq:harmonic} can be re-expressed as
\begin{equation*}
  h(z)=\frac{2}{\mu(B)} \Re\int_{B_0} \int_{B_1} f_{\alpha,\beta}\,
  d\mu(\beta)\, d\mu(\alpha),
\end{equation*}
and when $\mu(B_0)=\mu(B_1)$, this will be the real part of an
analytic function $f$ with $f(b)> 0$. 

Given $\alpha\in B_0$ and $\beta \in B_1$, let
\begin{equation}
  \label{eq:defpsi}
  \psi_{\alpha,\beta}=\frac{f_{\alpha,\beta}
    -f_{\alpha,\beta}(b)}{f_{\alpha,\beta}+f_{\alpha,\beta}(b)}.
\end{equation}
Then
\begin{equation*}
  f_{\alpha,\beta} = f_{\alpha,\beta}(b)
  \frac{\psi_{\alpha,\beta}+1}{\psi_{\alpha,\beta}-1}.
\end{equation*}
Note that $\psi_{\alpha,\beta}$ is unimodular on $B$, takes the
value $0$ at $b$, and in fact extends to an analytic function on a
region containing $\overline{\mathbb A}$.  Further,
$\psi_{\alpha,\beta}$ takes the value $1$ on $B$ precisely at those
points where $f_{\alpha,\beta}=\infty$; namely $\alpha$ and $\beta$.
Thus, $\psi_{\alpha,\beta}$ is two to one, and by the Maximum Modulus
Principle has two zeros in $\mathbb A$.  Since, the product of the
moduli of the zeros is $q$, the second zero is also on
the circle $\{|z|=\frac{q}{|b|}\}$.  The assumption that $b\neq
\sqrt{q}$ thus ensures that the zeros are distinct.

We claim that $\Theta^\prime=\{\psi_{\alpha,\beta}\}$ is a collection
of test functions for $H^\infty(\mathbb A)$; that is, that the unit
ball of $H^\infty(\mathbb A)$ is the same as the unit ball of
$H^\infty(\Theta^\prime)$.  One direction is nearly automatic. Since
$\Theta^\prime$ is a subset of the unit ball of $H^\infty(\mathbb A)$
it follows that the Szeg\H{o} kernel $s$ for $\mathbb A$ is in $\mathcal
K_{\mathcal K_{\Theta^\prime}}$. Thus, if $\varphi$ is in the unit ball of
$H^\infty(\mathcal K_{\Theta^\prime})$, then
\begin{equation*}
  \left((1-\varphi(x)\varphi(y)^*)s(x,y)\right) \geq 0.
\end{equation*}
Hence $\varphi$ is in the unit ball of $H^\infty(\mathbb A)$.  (In
general if $\Psi$ is contained in $\Psi^\prime$, then the unit ball of
$H^\infty(\mathcal K_{\Psi^\prime})$ is contained in the unit ball of
$H^\infty(\mathcal K_\Psi)$.)

 To prove the converse inclusion,
  suppose $\xi:\mathbb A\to
\mathbb D$ is analytic and $\xi(b)$ is real.  There exists a $\mu$ so
that
\begin{equation*}
  \begin{split}
    \frac{1+\xi}{1-\xi}
    =&\int_{B_0}\int_{B_1}  f_{\alpha,\beta}\, d\mu(\beta)\, d\mu(\alpha)\\
    =&\int_{B_0}\int_{B_1}
    f_{\alpha,\beta}(b)\frac{\psi_{\alpha,\beta}+1}{\psi_{\alpha,\beta}-1}
    \, d\mu(\beta)\, d\mu(\alpha).
  \end{split}
\end{equation*}
For $z,w \in \mathbb A$,
\begin{equation*}
  \begin{split}
    \frac{1+\xi}{1-\xi}(z)+\frac{1+\xi}{1-\xi}(w)^*
    &= 2 \frac{1}{1-\xi(z)} (1-\xi(z)\xi(w)^*)\frac{1}{1-\xi(w)^*}\\
    &=\int_{B_0}\int_{B_1} f_{\alpha,\beta}(b)
    \frac{1-\psi_{\alpha,\beta}(z)\psi_{\alpha,\beta}(w)^*}
    {(1-\psi_{\alpha,\beta}(z))(1-\psi_{\alpha,\beta}(w)^*)}
    \, d\mu(\beta)\, d\mu(\alpha).
  \end{split}
\end{equation*}
Thus, there exist functions $H_{\alpha,\beta}(z,w)$, analytic in $z$,
conjugate analytic in $w$ and continuous in $\alpha,\beta$ for fixed
$z,w$ so that
\begin{equation*}
  1-\xi(z)\xi(w)^* =\int_{B_0}\int_{B_1} H_{\alpha,\beta}(z,w)
  (1-\psi_{\alpha,\beta}(z)\psi_{\alpha,\beta}(w)^*) d\mu(\beta)\,
  d\mu(\alpha).
\end{equation*}
and the claim is proved.

There is some redundancy in our choice of test functions.  Given
$t\in\mathbb T$, let $\vartheta_t$ denote the function analytic in
$\mathbb A$, unimodular on $B$ with zeros at $b$ and $\frac{qt}{b}$
and with $\vartheta_t(1)=1$.  The collection
$\Theta=\{\vartheta_t:t\in\mathbb T\}$ is uniformly continuous in $t$
and $z$.  For each $\alpha,\beta$, there exist $t,\gamma$ in $\mathbb
T$ so that $\psi_{\alpha,\beta}=\gamma \vartheta_t$.  Thus, $\Theta$
is a totally bounded collection of test functions for
$H^\infty(\mathbb A)$.  (As an alternate, use the parameterization of
the unimodular functions with precisely two zeros in terms of theta
functions \cite{MR0335789}).

\begin{lemma}
  The collection $\Theta$ is a collection of test functions for
  $H^\infty(\mathbb A)$ and is compact in the norm topology of
  $H^\infty(\mathbb A)$.
\end{lemma}

Similar results for triply connected domains may be found in
\cite{MR2163865}.  See also a comment in \cite{MR1909298}.  The realization theorem,
Theorem~\ref{thm:main}, now reads as follows.

\begin{proposition}
  \label{prop:annulusalgebra}
  Suppose $\varphi:\mathbb A\to \mathbb C$.  The following are
  equivalent.
  \begin{enumerate}[(i)]
  \item $\varphi\in H^\infty(\mathbb A)$ with norm less than or
    equal to one;
  \item There is a positive kernel $\Gamma:\mathbb A\times \mathbb A
    \rightarrow C(\mathbb T)^*$ so that
    \begin{equation*}
      1-\varphi(z)\varphi(w)^*=\Gamma(z,w)(1-E(z)E(w)^*)
    \end{equation*}
    where $E(z)(\vartheta_t)=\vartheta_t(z)$; and
  \item there exists an auxiliary Hilbert space $\mathcal E$
    and an analytic function $\Phi:\mathbb A\to B(\mathcal E)$ whose
    values $\{\Phi(z)\}$ are commuting normal contraction operators
    and a unitary
    \begin{equation*}
      U=\begin{pmatrix} A &B\\ C&D \end{pmatrix} : 
      \begin{matrix} \mathcal E \\ \oplus \\ \mathbb C
      \end{matrix}\rightarrow 
      \begin{matrix} \mathcal E \\ \oplus \\ \mathbb C \end{matrix}
    \end{equation*} 
    so that $\varphi$ has the unitary colligation transfer function
    realization
    \begin{equation*}
      \varphi(z)=D+C\Phi(z)(I-A\Phi(z))^{-1}B.
    \end{equation*}     
  \end{enumerate}
\end{proposition}

We were able to reduce our original collection of test functions for
the annulus to $\Theta$.  It is reasonable to wonder if it is possible
to throw out even more.  The next proposition shows that the answer is
``no''.  It resembles a result of \cite{MR1386327} which says that in a sense
all of the Sarason/Abrahamse reproducing kernels for the annulus are
needed for Nevanlinna-Pick interpolation on the annulus.

\begin{proposition}
  \label{prop:needall}
  No proper closed subset of $\Theta$ is a set of test functions
  for $H^\infty(\mathbb A)$.
\end{proposition}

Note that a dense subset of the test functions will also be a set of
test functions, though in this case there is no real advantage to
taking such a set.  The situation will be quite different in the case
of the infinite polydisk, as we shall see.

\begin{proof}[Proof of Proposition~\ref{prop:needall}]
  Suppose $C$ is a proper subset of $\Theta$, and that $\vartheta_0 =
  \vartheta_{t_0}$ is not in $C$.

  Let $k:\mathbb A\times \mathbb A\to \mathbb C$ denote the Szeg\H{o}
  kernel for the annulus with respect to harmonic measure $\omega$ for
  the base point $b\in \mathbb A$ (recall that we assume $|b|\ne
  \sqrt{q}$).  Let $X$ denote the Schottky double of $\mathbb A$ and
  write $Jz$ for the twin of $z$ in the double.  According to Fay
  \cite{MR0335789} (see also \cite{MR1386327}), for each $a\in\mathbb A$ the kernel
  $k_a = k(\cdot,a)$ is meromorphic on $X$ with exactly two poles with
  the exception $k_b = 1$ (so in particular, $k(b,z) = k(z,b) = 1$ for
  all $z$).  Moreover, there exists a point $P$ in the complement of
  the closure of $\mathbb A$ so that $k_a$ has poles at $P$
  (independent of $a$) and $Ja$.

  The kernels
  \begin{equation*}
    \Delta_t(z,w)=(1-\vartheta_t(z)\vartheta_t(w)^*)k(z,w)
  \end{equation*}
  are positive and have rank two.  To see this, observe that $M_t$,
  the operator of multiplication by $\vartheta_t$ on $H^2(k)$, is an
  isometry, so $1-M_t M_t^*$ is the projection onto $\ker M_j^*$.
  Furthermore, if $b$ and $a_t$ are the two zeros of $\vartheta_t$
  (distinct since $|b|\ne \sqrt{q}$), then the identity $M_t^* k_w =
  \vartheta_t(w)^* k_w$ implies that $\ker M_t^* = \mathcal K_t =
  \mathrm{span}\,\{k_b, k_{a_t}\}$.  If we choose $f_t = k_b = 1$ and
  $g_t = k_{a_t} - k_b = k_{a_t} - 1$, then
  \begin{equation*}
    \Delta_t = f_tf_t^* + g_tg_t^* = 1 + g_tg_t^*.
  \end{equation*}
  It is useful to remark for later use that $g_t$ has the same poles
  as $k_{a_t}$; namely $P$ and $Ja_t$.

  Choose three distinct points $z_1,z_2,z_3\in\mathbb A$ and consider
  the Agler-Pick interpolation problem of finding a $\varphi\in
  H^\infty(\mathbb A)$ so that $\varphi(z_t)=\vartheta_0(z_t)$.  The fact
  that the $3\times 3$ matrix
  \begin{equation*}
    \begin{pmatrix} k(z_\ell,z_m)(1-\vartheta_0(z_\ell)\vartheta_0(z_m)^*)
    \end{pmatrix}_{j,m=1}^3
  \end{equation*}
  has rank two implies that this interpolation problem has a unique
  solution, namely $\varphi=\vartheta_0$.  On the other hand, there is a
  bounded positive measure $\mu$ on $C$ so that $\varphi$ has a
  realization of the form in Theorem \ref{thm:APint} with $n=3$.  The
  usual computations convert that realization to
  \begin{equation*}
    1-\vartheta_0(z)\vartheta_0(w)^* 
    = \int_C \sum_{\nu=1}^3 h_\nu(z,\vartheta)
    h_\nu(w,\vartheta)^*(1-\vartheta(z)\vartheta(w)^*)\, d\mu(\vartheta)
  \end{equation*} 
  for functions $h_t(z,\cdot)\in L^2(\mu)$.  In particular,
  multiplying through by $k(z,w)$ gives
  \begin{equation*}
    \Delta_0(z,w) = \int_C \sum_\nu h_\nu(z,\vartheta) \Delta_\vartheta(z,w)
    h_\nu(w,\vartheta)^* \, d\mu(\vartheta)
  \end{equation*}
  where $\Delta_\vartheta(z,w) = (1-\vartheta(z)\vartheta(w)^*)k(z,w)$.

  Fix $z$.  Since $\Delta_\vartheta(z,z) \geq 0$ $\mu$ a.s., given
  $\delta > 0$, there is a set $C'\subset C$ and a constant
  $c_\delta>0$ such that $\mu(C-C') < \delta$ and for all $z\in\mathbb
  A$,
  \begin{equation*}
    \Delta_0(z,z) \geq c_\delta h_\nu(z,\vartheta) \Delta_\vartheta(z,w)
    h_\nu(w,\vartheta)^*, \qquad \vartheta\in C'.
  \end{equation*}
  Then using the factorization of the $\Delta$'s given above, by
  Douglas' lemma there are constants $c_k$, $k=1,2,3,4$, such that for
  fixed $\vartheta\in C'$,
  \begin{equation*}
    \begin{split}
      h_\nu(\cdot,\vartheta) &= c_1 + c_2 g_{a_0} \\
      h_\nu(\cdot,\vartheta)g_{a_t} &= c_3 + c_4 g_{a_0}.
    \end{split}
  \end{equation*}
  Since the kernels extend meromorphically to $X$, the same is true
  for $h_\nu(\cdot,\vartheta)$ by the first equation.  That equation
  also implies that either $h_\nu(\cdot,\vartheta)$ is constant or
  that it has the same poles as $g_{a_0}$; that is, simple poles at
  $P$ and $a_0$.  If $h_\nu(\cdot,\vartheta)$ is not constant, then
  the left side of the second equation has a double pole at $P$, while
  the right only has a single pole.  Hence $h_\nu(\cdot,\vartheta)$
  must be constant.  If it is a nonzero constant, the second equation
  would imply that the poles of $g_{a_t}$ and $g_{a_0}$ agree, and in
  particular, that $a_t = a_0$, contradicting the assumption that
  $\vartheta_0 \notin C$.  Hence $h_\nu(\cdot,\vartheta) = 0$.

  Taking $\delta$ going to $0$, we see that the subset of $C'$ on
  which $h_\nu(\cdot,\vartheta)$ is nonzero has $\mu$ measure zero,
  yielding a contradiction.
\end{proof}

\subsection{The infinite polydisk}
\label{subsec:F}

By the infinite polydisk ${\mathbb D}^\infty$, 
we mean the open unit ball of $C_b(\mathbb N)$. Thus,
\begin{equation*}
  {\mathbb D}^\infty =\{z:\mathbb N\to \mathbb D: \sup
  \{|z(n)|:n\}<1\}. 
\end{equation*}

Let $e_n$ denote the function $e_n :{\mathbb D}^\infty \to \mathbb C$
given by $e_n(z)=z(n)$. The set of test functions
$\Psi=\{e_n:n\in\mathbb N\}$ is topologized by the inclusion
$\Psi\subset B({\mathbb D}^\infty, \overline{\mathbb D})$.  The spaces
$\Psi$ and $\mathbb N$ are homeomorphic and hence $\beta\Psi$ is
identified with $\beta \mathbb N$.

A $\chi\in \beta\mathbb N\backslash \mathbb N$ determines a function
$\varphi_\chi :X \to \mathbb D$ given by
\begin{equation*}
  \varphi_\chi (z) = z(\chi),
\end{equation*}
where we have identified $z\in X$ with its unique extension to a
continuous function $z:\beta\mathbb N \to \mathbb C$.  This
identification follows from the general discussion of $\Psi$ and
$\beta\Psi$ in subsection~\ref{subsec:test-funct-eval}.  Further,
$\varphi_\chi$ is in the unit ball of $H^\infty(\Psi)$.

Theorem~\ref{thm:main} now implies that there is a positive kernel
with entries in $C(\overline{{\mathbb D}^\infty})^*$ such that
\begin{equation*}
  1-\varphi(z)\varphi(w)^* = \Gamma(z,w) (1-E(z)E(w)^*)\geq 0.
\end{equation*}
In this case there is a clear choice for $\Gamma$; namely,
$\Gamma(z,w) = \gamma$, where $\gamma(e) = e(\chi)$ for $e\in
C(\overline{{\mathbb D}}^{\mathbb N})$.

\begin{proposition}
 \label{prop:infinite-polydisc-no-no}
 There \textbf{do not} exist positive kernels $\Gamma_n :X\times X\to
 \mathbb C$ such that
 \begin{equation*}
   1-\varphi_\chi(z)\varphi_\chi(w)^* =\sum_n \Gamma_n(z,w)
   (1-z(n)w(n)^*).
 \end{equation*}
 Similarly, if $C$ is any closed subset of $\beta\mathbb N$ with
 $\chi\notin C$, then there \textbf{does not} exist a positive
 $\Gamma:X\times X\to C(C)^*$ such that
 \begin{equation*}
   1-\varphi_\chi(z)\varphi_\chi(w)^* =\Gamma(z,w)(I-E(z)E(w)^*).
 \end{equation*}
\end{proposition}

\begin{proof}
  For the first part observe that if $z(n)$ converges to $L$ as
  $n\to\infty$, then $\chi(z)=L$ for any $\chi\in \beta\mathbb
  N\backslash\mathbb N$.

  Choose $z(n)=\sqrt{\frac{1}{2}\left(1-\frac{1}{n+1}\right)}$ we have
  $\frac{1}{2} = z(\chi) > z(n)\ge 0$ for all $n$.  Let $0$ denote the
  zero sequence.  Suppose that the first representation in the
  proposition holds.  Then
  \begin{equation*}
    \frac{1}{2} = \sum_n \Gamma_n(z,z)\tfrac{1}{2}\left(1 +
      \tfrac{1}{n+1}\right) > \tfrac{1}{2}\sum_n \Gamma_n(z,z),
  \end{equation*}
  and so $\sum_n \Gamma_n(z,z) < 1$.  Obviously $\sum_n \Gamma_n(z,0)
  = \sum_n \Gamma_n(0,z) = 1$.  Also, for each $n$,
  \begin{equation*}
    \begin{pmatrix}  \Gamma_n(z,z) & \Gamma_n(z,0) \\ 
      \Gamma_n(0,z) & \Gamma_n(0,0) \end{pmatrix} \geq 0,
  \end{equation*}
  so
  \begin{equation*}
    \sum_n \begin{pmatrix}  \Gamma_n(z,z) & \Gamma_n(z,0) \\ 
      \Gamma_n(0,z) & \Gamma_n(0,0) \end{pmatrix} = 
    \begin{pmatrix}  \sum_n \Gamma_n(z,z) & 1 \\ 
      1 & 1 \end{pmatrix}\geq 0,
  \end{equation*}
  and thus $\sum_n \Gamma_n(z,z) \geq 1$, a contradiction.

  The second part of the Proposition is proved similarly, in this
  situation choosing a function $z$ with $z(C)=0$,
  $z(\chi)=\sqrt{\tfrac{1}{2}}$, and $0\le z \le
  \sqrt{\tfrac{1}{2}}$ (such a function exists since $\beta\mathbb N$
  is Tychonov).
\end{proof}

By the way, if we define $P_n$ as the projection of $f\in {\mathbb
  D}^{\mathbb N}$ onto its first $n$ components, then despite the fact
that $P_n f$ converges pointwise to $f$, $\varphi(P_n f) = 0$, and so
obviously does not converge to $\varphi(f)$ in general.  However this
does not contradict Proposition~\ref{prop:pw-simple}, since this only
says that there is some net of simple representations converging to
$\varphi$, and obviously this is not one!

\subsection{Further Examples}
\label{subsec:further-examples}

We end with a few examples illustrating some of the pitfalls into
which an unwary applicant of the results presented can fall.

The following example shows that it is sometimes important and natural
to use the compactification of $\Psi$, and illustrates more simply the
phenomena observed with the infinite polydisk.

\subsubsection{Example 1.}
\label{subsubsec:example1}

Choose $X$ equal the unit disk $\mathbb D$ and let, for $n=1,2,\dots,$
$\psi_n(z)=\left(\sqrt{1-\frac{1}{n}}\right)z$.  The collection $\Psi
= \{\psi_n:n\}$ is a set of test functions for $H^\infty(\mathbb D)$
and the function $\xi(z)=z$ is in $H^\infty(\mathbb D)$ with
$\|\xi\|=1$.

\begin{lemma} There do not exist positive kernels $\Gamma_n$ so that
  \begin{equation}
    \label{eq:noPsD}
    1-zw^*=\sum_{n\in\mathbb N} \Gamma_n(z,w)(1-\psi_n(z)\psi_n(w)^*).
  \end{equation}
\end{lemma}

\begin{proof}
  Suppose (\ref{eq:noPsD}) holds for some positive kernels $\Gamma_j$.
  Divide through by $1-zw^*$ to obtain,
\begin{equation*}
  1=\sum \Gamma_n(z,w) +\frac{1}{n}\frac{\Gamma_n(z,w)}{1-zw^*}.
\end{equation*}
Note that the left side is the rank one positive matrix $[1]$
consisting of all $1$'s, and also each term on the right side is
positive.  Hence each term on the right side is a nonnegative constant
multiple of $[1]$, which is clearly a contradiction.
\end{proof}

Interestingly, if we had used the function
$\left(\sqrt{1-\frac{1}{n}}\right)z$ instead of $z$, then there is an
obvious choice for the $\Gamma_k$'s; namely, $\Gamma_n = [1]$ and all
others equal to $0$.  Furthermore, for a finite set $F \subset X$, the
matrices
\begin{equation*}
  \left(1-\left(1-\tfrac{1}{n}\right)zw^*\right)_{z,w\in F}
\end{equation*}
converge to
\begin{equation*}
  \left(1-zw^*\right)_{z,w\in F},
\end{equation*}
and so by the proof of Lemma~\ref{lem:closedcone}, there must be a
$\Gamma$ such that $1-zw = \Gamma(z,w)(1-E(z)E(w)^*)$ (the proof
of the lemma makes no use, either explicit or implicit, of the
compactness of the set of test functions).

The kernel $\Gamma$ has entries which are continuous functions over
$C_b(\Psi)$, and this includes the point evaluations which we tried to
use above.  However there are point evaluations we have not considered
--- the ones coming from points in the Stone-\v{C}ech compactification
$\beta\Psi$ of $\Psi$.  In this case, this agrees with the one point
compactification where we add the function $\psi_\infty(z) = z$.  The
functions $E(z)$ extend uniquely to $\Psi$, and if the positive linear
functional $\gamma\in C_b(\Psi_0)^*$ is defined by $\gamma(e) =
e(\psi_\infty)$, then the choice $\Gamma(x,y) = \gamma$ for all $x,y$
gets us out of our quandary.

In this example, it is clear that there would be no harm (in fact, it
would be to our advantage) to include the functions in the
Stone-\v{C}ech compactification of the set of test functions,
particularly since it does not much effect the size of the set of test
functions.  This is in stark contrast to the case of the infinite
polydisk, where the compactification increases the set size from being
countable to at least having cardinality of $2^{\mathfrak c}$.

\subsubsection{Example 2.}

Let $X=\{x_1,x_2\}$ denote a two point set and define $\psi(x_1)=0$
and $\psi(x_2)=1$.  The set $\Psi=\{\psi\}$ is a set of test functions
for $C(X)$.  However, the function $\tilde{\psi}=1-\psi$ is in the
unit ball of $C(X)$, but there does not exist a positive $\Gamma$ such
that
\begin{equation}
  \label{eq:noGamma3}
  1-\tilde{\psi}(x)\tilde{\psi}(y)^* = \Gamma(x,y)(1-\psi(x)\psi(y)^*).
\end{equation}
The remainder of this subsection is devoted to these assertions.

Suppose $k\in\mathcal K_\Psi$; that is, $k$ is a positive kernel and
the $2\times 2$ matrix
\begin{equation*}
    \begin{pmatrix} (1-\psi(x_1)\psi(x_1)^*)k(x_1,x_1) &
      (1-\psi(x_1)\psi(x_2)^*)k(x_1,x_2)\\
      (1-\psi(x_2)\psi(x_1)^*)k(x_2,x_1) &
      (1-\psi(x_2)\psi(x_2)^*)k(x_2,x_2)
    \end{pmatrix} =
    \begin{pmatrix} k(x_1,x_1) & k(x_1,x_2)\\ k(x_2,x_1) &
      0\end{pmatrix}
\end{equation*}
is positive.  It follows that $k(x_1,x_1)$ and $k(x_2,x_2)$ are
nonnegative and $k(x_1,x_2)=0=k(x_2,x_1)$.  Now for
$\varphi:X\to\mathbb C$ it is readily verified that
$(1-\varphi(x)\varphi(y)^*)k(x,y)$ is positive for all such $k$ if and
only if $|\varphi(x_j)|\le 1$ for $j=1,2$.  Hence $\HP = C(X)$ and
$\tilde{\psi}$ is contractive.

Since the $(x_2,x_2)$ entry of the left hand side of equation
(\ref{eq:noGamma3}) is $1$, but the same entry on the right hand side
of this equation is $0$, no such $\Gamma$ exists.  So what went wrong?
To begin with, $\psi$ violates condition (i) for a test function.
However this is not so serious in this case.  More to the point, it
also violates (ii), since if we choose the set $F = \{x_2\}$, $\Psi|F$
does not generate $P(F)$.  We could fix this either by taking a
quotient or by adding another test function.  A natural choice is
$\tilde\psi$; sadly this still violates condition (i).

\bibliographystyle{plain}
\bibliography{test_functions}

\end{document}